\documentclass[draft]{amsart}
\usepackage{amssymb,latexsym}
\usepackage[mathscr]{eucal}
\usepackage{graphics}
\usepackage{ulem} 
\usepackage{verbatim}

\theoremstyle{plain}
\newtheorem{theorem}{Theorem}[section]
\newtheorem{corollary}[theorem]{Corollary}
\newtheorem{lemma}[theorem]{Lemma}
\newtheorem{proposition}[theorem]{Proposition}
\newtheorem{remark}[theorem]{Remark}
\newtheorem{definition}[theorem]{Definition}
\newtheorem{example}[theorem]{Example}

\def\sideremark#1{\ifvmode\leavevmode\fi\vadjust{\vbox to0pt{\vss
\hbox to 0pt{\hskip\hsize\hskip1em
\vbox{\hsize2cm\tiny\raggedright\pretolerance10000
\noindent#1\hfill}\hss}\vbox to8pt{\vfil}\vss}}}

\newcommand{\be}{\begin{equation}\label}
\newcommand{\ee}{\end{equation}}
\newcommand{\bq}{\begin{equation*}}
\newcommand{\eq}{\end{equation*}}
\newcommand{\ba}{\begin{align*}}
\newcommand{\ea}{\end{align*}}
\newcommand{\bp}{\begin{proof}}
\newcommand{\ep}{\end{proof}}
\newcommand{\bL}{\begin{lemma}\label}
\newcommand{\eL}{\end{lemma}}
\newcommand{\bP}{\begin{proposition}\label}
\newcommand{\eP}{\end{proposition}}
\newcommand{\bC}{\begin{corollary}\label}
\newcommand{\eC}{\end{corollary}}
\newcommand{\bT}{\begin{theorem}\label}
\newcommand{\eT}{\end{theorem}}
\newcommand{\bR}{\begin{remark}\label}
\newcommand{\eR}{\end{remark}}
\newcommand{\bD}{\begin{definition}\label}
\newcommand{\eD}{\end{definition}}
\newcommand{\bE}{\begin{example}\label}
\newcommand{\eE}{\end{example}}

\DeclareMathOperator{\tr}{Tr}
\DeclareMathOperator{\diag}{diag}
\DeclareMathOperator{\Adm}{Adm}

\date{10/13/2017}
\thanks{This work was partially supported by the Simons Foundation (grant No 245660 to Victor Kaftal)}
\author{Victor Kaftal}

\address{Department of Mathematics\\
University of Cincinnati\\
P. O. Box 210025\\
Cincinnati, OH\\
45221-0025\\
USA}

\email{victor.kaftal@uc.edu}

\author{David R. Larson}

\address{Texas A\&M University\\
Department of Mathematics\\
College Station, TX 77843\\
USA}

\email{larson@math.tamu.edu}


\begin{document}

\title{Admissible sequences for positive operators}
\begin{abstract}
A scalar sequence $\xi$ is said to be admissible for a positive operator $A$ if $A= \sum\xi_jP_j$ for some rank-one projections $P_j$, or, equivalently, if $\xi$ is the diagonal of $VAV^*$ for some partial isometry $V$ having as domain the closure of the range of $A$. The main result of this paper is that if $\xi$ is a non summable sequence in $[0,1]$ that satisfies the Kadison condition  that either  $\sum\{\xi_i \mid \xi_i\le \frac{1}{2}\}+ \sum\{(1-\xi_i) \mid \xi_i>   \frac{1}{2}\} = \infty$ or $\sum\{\xi_i \mid \xi_i\le \frac{1}{2}\}- \sum\{(1-\xi_i) \mid \xi_i>   \frac{1}{2}\} \in \mathbb Z$ and if $A$ is the sum of infinitely many projections (converging in the SOT) then $\xi$ is admissible for $A$. This result extends Kadison's carpenter's theorem and provides an independent proof of it.
\end{abstract}
\maketitle
\section{Introduction}
The study of diagonals of operators started when Schur  (\cite{Si23}) proved that if $\xi$ is the diagonal sequence of a (finite) selfadjoint matrix and $\lambda(A) $ is the eigenvalue list with multiplicity  of $A$, then $\xi$ is majorized by $\lambda(A)$ ($\xi\prec \lambda(A)$) and then  Horn (\cite{Ha54}) proved that this condition is also sufficient. Although majorization theory has played an important role in various areas of analysis and algebra since the early 1900's, there is no universally accepted notation for it, so we need to state explicitly our notations:

\bD{D:major}
Let $\xi,\lambda\in \mathbb R^n$ and $\xi^*,\lambda^*$ denote their monotone
non-increasing rearrangements. Then 
$\lambda$ majorizes $\xi$
($\xi\prec \lambda$) if
$\sum_{j=1}^k\xi^*_j\le \sum_{j=1}^k\lambda^*_j$ for every $1\le k\le
n$ 
and  $\sum_{j=1}^n\xi_j= \sum_{j=1}^n\lambda_j.$ \\
Let $\xi, \lambda\in (c_o)_+$  and $\xi^*,\lambda^*$ denote their monotone
non-increasing rearrangements. Then $\xi\prec \lambda$ if $\sum_{j=1}^k\xi^*_j\le \sum_{j=1}^k\lambda^*_j$ for every $k\ge 1$
and  $\sum_{j=1}^\infty \xi_j= \sum_{j=1}^\infty \lambda_j$.
\eD
For convenience, if  $\xi\in \mathbb R_+^n$ and $ \lambda\in \mathbb R_+^m$ then we will still say that $\xi\prec \lambda$ if $\tilde \xi\prec \tilde \lambda$ where $ \tilde\xi, \tilde \lambda \in \mathbb R_+^N$ for $N= \max\{m,n\}$  are the sequences obtained by completing $\xi$ or $\lambda$ with zeros. 

The Schur--Horn theorem was extended to infinite dimensional positive compact operators in \cite {aM64}, \cite {GiMa64}, \cite {AK02}, \cite {KaftalWeiss2010}, \cite { LjWg15}. More precisely, given  a separable complex  Hilbert space $H$ with an othonormal basis $\{e_n\}$, and the conditional expectation $E$ on the masa of operators diagonal with respect to that basis (i.e., the operation of ``taking the main diagonal"), and given an operator $A\in K(H)_+$, the majorization condition $\xi\prec\lambda(A)$ characterizes the diagonals of the {\it partial isometry orbit} of $A$, $\mathscr V(A)= \{ VAV^*\mid V^*V= R_A\}$ where $R_A$ denotes the range projection of $A$. The characterization of the diagonals of the {\it unitary orbit} $\mathscr U(A)$ of a positive compact operator is more delicate and its study is not yet complete (e.g., \cite {KaftalWeiss2010}, \cite {LjWg15}.)

Central to the topic of our paper is the seminal work of Kadison (\cite{Kr02a}, \cite{Kr02b}) characterizing the diagonals of projections in terms of the following property:   

\bD {D: Kad cond} A sequence $\xi:=\langle\xi_i\rangle $  (finite or infinite)  is said to satisfy the Kadison condition if $0\le \xi_i\le 1$ for every $i$ and  either 
\ba \sum\{\xi_i \mid \xi_i\le \frac{1}{2}\}+ \sum\{(1-\xi_i) \mid \xi_i>   \frac{1}{2}\} &= \infty\\\intertext{or}\sum\{\xi_i \mid \xi_i\le \frac{1}{2}\}- \sum\{(1-\xi_i) \mid \xi_i>   \frac{1}{2}\}& \in \mathbb Z.\end{align*}
\eD
\bR{R:1/2} As it has been often observed, a sequence $\xi$ satisfies the Kadison condition if and only if it satisfies the condition where $ \frac{1}{2}$ is replaced by a number $0< \alpha< 1$ (but the integer may change.) Similarly, the conditions $\le   \frac{1}{2}$ and  $> \frac{1}{2}$ can be replaced respectively by the conditions $<  \frac{1}{2}$ and $\ge \frac{1}{2}$. \eR

\bT{T:Kadison}\cite[Theorem 15]{Kr02b}
A sequence $\xi$ satisfies the Kadison condition if and only if $\diag \xi= E(R)$ for some projection $R$.
\eT

The necessity part of the theorem has seen further work in \cite{Am15},\cite {AW07}, \cite {KvLj08}) but we will say no more about it as it is less related to our present paper.

The sufficiency part of the theorem is commonly called the carpenter's theorem.
Its proof is quite hard although it has seen some simplifications in \cite{Am15} and \cite{BJCarpenter}. We find it convenient to reformulate Kadison's carpenter's  theorem as follows:  
\bT{T:carpenter} If $\xi$ satisfies the  Kadison condition, and $P$ is a projection such that $\sum_j\xi = \tr(P)$  then $\diag \xi\in E(\mathscr V(P))$. 
 \eT
\bp By Theorem \ref {T:Kadison}, $\diag \xi= E(Q)$ for some  projection $Q$. But then $\tr(P)=\sum_j\xi = \tr(Q)$ and hence  $Q\sim P$, i.e., $Q\in \mathscr V(P)$.\ep

Of course, if $\xi$ satisfies the Kadison condition, then $\sum_j\xi = \tr(P)$ is also necessary for $\diag \xi\in E(\mathscr V(P))$. 
Kadison's characterization of the diagonals of projections, or, more generally, of selfadjoint operators with two point spectrum, was then extended to a characterization of the diagonals of selfadjoint operators with finite spectrum by Bownik and Jasper (\cite {BJSchHor}, \cite{BJdiag}).  Arveson extended the necessity of  the condition to normal operators with finite spectrum (\cite {AW07}).

\

An apparently unrelated line of research was started in \cite {Fp69} by Fillmore who characterized  the positive finite rank operators on a separable Hilbert space that are sums of projections. 
More recent work on frames with prescribed vector norms can be seen as an extension of Fillmore's result (\cite {CFKLT06}, \cite{CL02}, \cite{CL10}, and \cite {KkLd04}). We quote the following result that is important for our paper: 
\bT {T:KkLd}  \cite [Theorem 2] {KkLd04}
Let $A$ be positive operator with rank $r$, let $\xi  \in (\mathbb R^n)_+$ for some $n\ge r$. Then  $A= \sum_{j=1}^n\xi_jP_j$ for some  rank-one projections $P_j$, if and only if $\xi\prec \lambda(A)$.
\eT

Notice that the same majorization condition appears both here and in the (finite dimensional) Schur--Horn theorem. It may therefore be interesting to further examine the direct link between diagonals of a positive operator and  decompositions of that operator into positive linear combinations of projections. This link holds also in infinite dimensions and has been presented implicitly  in \cite {AMRS} and explicitly in \cite [Proposition 3.1] {KNZ strong sums}, in an early draft of this paper, in \cite [Proposition 5.1]{KNZFiniteSumsVNA},  and in \cite [Proposition 3.5] {DraganKaftal2016}. Given its simplicity, usefulness, and its central role in this paper, we present it here as well. In order to be able to handle both finite and infinite sequences, we embed $\mathbb R^n$ into $\ell^\infty.$
\bP{P: bridge} 
Let 
$A\in B(H)_+$  and $\xi\in \ell^\infty_+$. Then the following conditions are equivalent 
\item[(i)] There is a partial isometry $V$ with $V^*V = R_A$ for which $\diag \xi = E(VAV^*)$.
\item [(ii)] There is a sequence of rank one projections $P_j$ such that $A= \sum_{j=1}^\infty \xi_jP_j$, where the series converges in the strong operator topology.
\eP
The proof of this proposition will be presented in the next section for the convenience of the reader.  
Let us mention a few cases that illustrate the potential usefulness of Proposition \ref {P: bridge}. By using it, one can obtain Theorem \ref {T:KkLd}
directly from the celebrated Schur--Horn theorem. Conversely, Theorem \ref {T:KkLd} provides a proof of the Schur--Horn theorem which is noticeably simpler than the original proof by Horn \cite {Ha54}, and in our view, perhaps simpler also than some of the other more recent proofs of the Schur--Horn theorem. 

Another simple consequence of Proposition \ref {P: bridge} is that an operator $A$ is the sum of projections (strongly converging, if infinitely many),  if and only if $I\in E(\mathscr V(A))$. In fact, a  simple modification of the proof of Proposition \ref {P: bridge} (see \cite [Proposition 5.1] {KNZFiniteSumsVNA}) shows that $A$ is the sum of $n$ projections  if and only if there is an $n\times n$ matrix decomposition of $A$ with the block diagonal being the identity and this fact provides another proof of \cite[Theorem 1.2] {ChoiWu2014}.

Sums of projections  are  the frame operators of  unit norm frames or more precisely of Bessel sequences. Their study in this context was started in  \cite {DFKLOW} where it was shown that a sufficient condition is that $\|A\|_{e}> 1$ (essential norm). Further advances in that direction were obtained in \cite {AMRS} and the full characterization  was  obtained in \cite{KNZ strong sums}:
\bT{T;sums of proj}\cite [Theorem 1.1] {KNZ strong sums} Let $A\in B(H)_+$. Then $A$ is the sum of projections if and only if  either $\tr((A-I)_+)=\infty$ or $\tr((A-I)_+)- \tr((I-A)_-R_A)\in \mathbb N$ 
\eT
Here and in the following, by  $\mathbb N$ we mean $\{0,1,2, \cdots\}$. 
In this paper we will focus  on the sequences associated to a positive operator $A$ by the (equivalent) properties of Proposition \ref {P: bridge}. 
\bD{D:Adm} A sequence $\xi \in\ell^\infty_+$ is said to be {\it admissible} for an operator $A\in B(H)_+$ if $A= \sum_{j} \xi_jP_j$ for some sequence of rank-one projections, where the series converges in the strong operator topology if infinitely many. We denote by $\Adm(A)$ the collection of all admissible sequences for $A$.
\eD

Reformulated in terms of this notation, the Schur--Horn theorem states that if $A$ is a finite matrix then $\Adm(A)=\{\xi\in \ell^\infty_+\mid \xi\prec \lambda(A)\}$ and precisely the same characterization holds for every $A\in K(H)_+$ (\cite [Proposition 6.4] {KaftalWeiss2010}). All infinite projections share the same admissible sequences and Kadison's Theorem \ref {T:Kadison} states that these are precisely the nonsummable sequences satisfying the Kadison condition. 

The question that initiated our project a few years ago  was whether  Proposition \ref {P: bridge} could provide an alternative approach to the proof of carpenter's Theorem \ref {T:carpenter}. We were successful and, more interestingly, we found that the same techniques could be applied to sums of projections. The following is the main theorem of our paper.

\bT{T:main theorem}
Let $\xi$ be a sequence satisfying the Kadison condition and let $A$ be a sum of projections such that  $\tr(A)=\sum_j \xi_j$. Then $\xi\in \Adm(A)$.
\eT

Our paper is structured as follows.  In Section 2 we prove Proposition \ref {P: bridge} and then derive some preliminary results on admissible sequences and a corollary of the main theorem illustrating how rich is the admissible class of a noncompact positive operator.

In Section 3 we present the proof of the main theorem which is split into a number of cases. In most instances the proof is obtained by splitting an infinite sum of projections $A=\sum_jE_j$ into a sum $A=\sum A_j$ of finite rank operators
$$A_j= (1-r_j)E_{n_j} + E_{n_j+1}+\cdots+  E_{n_{j+1}-1}+ r_{j+1}E_{n_{j+1}}$$
 chosen so that 
\be {e:major}\langle \xi_{m_j},\cdots,  \xi_{m_{j+1}}\rangle\prec \langle 1-r_j, 1, \cdots, 1,r_{j+1}\rangle\ee
 and hence by Theorem \ref {T:KkLd} or more precisely, by its Corollary \ref {C:from KornLars}, $A_j= \xi_{m_j}P_{m_j}+ \cdots  \xi_{m_{j+1}}P_{m_{j+1}}$, where $P_i$ are rank one projections. The technical difficulty is in the appropriate choice of the sequences $n_j$, $m_j$, and the remainders $r_j$ so to guarantee the majorization in (\ref {e:major}).
One notable exception  is the case when $\xi$ clusters summably under 1 (Lemma \ref  {L:M finite N infinite}), which requires a completely different technique in the decomposition and especially in proving convergence. Notice that all the proofs do not depend on Kadison's Theorem \ref {T:carpenter}  and thus provide  a new proof of it. However, it may be interesting to notice that if $A$ is itself a projection, both Lemma \ref  {L:M finite N infinite} as well as Lemma \ref {L: infinite sum of lambda} are immediate consequences  of the previously obtained (and simpler) lemmas applied to $I-A$. This greatly simplifies the proof for the case of a projection.

The first named author wishes to thank John Jasper for useful discussions on this paper during his stay at the University of Cincinnati and thank also Jireh Loreaux for many useful comments.

 \section{Preliminary results on admissible sequences.}

The  motivation for  the notion of admissible sequences arises from Proposition \ref {P: bridge}. We will present therefore a short proof of it.
\bp [Proof of Proposition \ref {P: bridge}]
(i)$\Rightarrow$(ii). By the condition $\diag \xi = E(VAV^*)$, $$\|A^{1/2}V^*e_j\|^2= (VAV^*e_j,e_j)=\xi_j\quad \forall \,  j.$$
Set $\Lambda:=\{j\in \mathbb N\mid \xi_j\ne 0\}$ and for $j\in \Lambda$, let $v_j:= \frac{1}{\sqrt {\xi_j}} A^{1/2}V^*e_j$. Then $\|v_j\|=1$ and $P_j:=v_j\otimes v_j$  is a rank one projection. Moreover,
\ba A&=A^{1/2}V^*VA^{1/2}= A^{1/2}V^*\Big(\sum_{j=1}^\infty e_j\otimes e_j\Big)VA^{1/2}\\
&=  \sum_{j=1}^\infty A^{1/2}V^*e_j\otimes A^{1/2}V^*e_j= \sum_{j\in \Lambda}A^{1/2}V^*e_j\otimes A^{1/2}V^*e_j\\
&=  \sum_{j\in \Lambda}\xi_jP_j
\end{align*}
where the strong convergence of the series $\sum_{j\in \Lambda}\xi_jP_j$ is an immediate consequence of the strong convergence of the series $\sum_{j=1}^\infty e_j\otimes e_j=1$.\\
(ii)$\Rightarrow$(i) Let $v_j$ be a unit vector in the range of the projection $P_j$. By hypothesis $A=\sum_{j=1}^\infty \sqrt{\xi_j}v_j\otimes \sqrt{\xi_j}v_j$ where the series converges strongly. Then it is routine to see that the series $B:= \sum_{j=1}^\infty \sqrt{\xi_j}e_j\otimes v_j$ also converges strongly and that $B^*e_j= \sqrt{\xi_j}v_j$. Moreover,
$$ B^*B= \sum_{i,j=1}^\infty (\sqrt{\xi_i}v_i\otimes e_i )(\sqrt{\xi_j}e_j\otimes v_j)=  \sum_{j=1}^\infty \xi_jv_j\otimes v_j =A.$$
Let $B=VA^{1/2}$ be the polar decomposition of $B$. Then $V^*V=R_A$ and for all $j$,
$$(VAV^*e_j,e_j)= \|B^*e_j\|^2= \| \sqrt{\xi_j}v_j\|^2= \xi_j.$$
Thus $E(VAV^*)=\diag \xi$.
\ep
\bR{R:frames} 
The objects that enter this proof are familiar in frame theory. Indeed $\{\sqrt{\xi_j}v_j\}$ is a Bessel sequence (not necessarily a frame since  we do not need to assume that  $A$ is invertible), $A$ is its frame operator,  $VAV^*$ is its Gram matrix, $B$ is its analysis operator and  $V$ is the analysis operator of the associated Parseval frame.
\eR

\bR{R:other equivalence}
The equivalent conditions in Proposition \ref {P: bridge} are also equivalent to\\
(iii) $\diag \xi \oplus 0= E(U(A\oplus 0)U^*)$ for some unitary $U$ and some direct summand $0$.\\
The equivalence of (iii) with the condition (ii) Proposition \ref {P: bridge}  has been obtained in \cite [Proposition 4.5] {AMRS} for the case of that $A$ is a  frame operator, but the proof did not depend on the invertibility of $A$.
\eR

Clearly, due to the unconditional convergences of the series $A= \sum_{j=1}^\infty \xi_jP_j$, the admissibility of a sequence $\xi$ is not affected by permuting the sequence nor by adding or deleting zeros from the sequence. It is thus convenient to embed $\mathbb R^n$ into $\ell^\infty$ so to handle at the same time  both  finite and infinite sequences.  It is also convenient to use the following notation:
$$\text{If }\xi, \eta \in \ell^\infty, \text{ then }\xi\oplus \eta= <\xi_1, \eta_1, \xi_2, \eta_2, \dots>.$$ with the obvious modification if one or both the sequences $\xi$ or $ \eta$ are finite.

We collect here some elementary properties of the classes of admissible sequences.
\bL{L: trivial properties} Let $A, B\in B(H)_+$. Then
\begin{enumerate}
\item [(i)]  $\Adm(A)\subset [0, ||A||]^{\mathbb N}$;
\item [(ii)] $\Adm(A)= \Adm(VAV^*)$ for every partial isometry  $V$  with  $V^*V \ge R_A$;
\item [(iii)] $\Pi \Adm(A) =\Adm(A)$ for every permutation $\Pi$;
\item [(iv)] $\Adm(A) = \Adm(A\oplus 0) = \Adm(A)\oplus 0$ for $0$ operator summands and $0$ sequence summands of any size;
\item [(v)]$\gamma \Adm(A)= \Adm(\gamma A)$ for every $\gamma\ge 0$. In particular, $ \Adm(0)= \{0\}$;
\item [(vi)] $\Adm(A)\oplus \Adm(B)\subset \Adm(A+B)$ and the inclusions can be proper;
\item [(vii)] If $\xi \in \Adm(A)$, then $\sum_{j=1}^\infty \xi_j= \tr(A)$;
\item [(viii)] If $\xi, \eta  \in \Adm(A)$ and $0\le t\le 1$, then $t\xi\oplus(1-t)\eta \in \Adm(A)$;
\item [(ix)] If $\xi\in \Adm(A)$ and $0\le \eta\le \xi$ then $\eta\oplus (\xi-\eta)\in \Adm(A).$
\end{enumerate}
\eL
To illustrate how the above facts follow easily from Definition \ref {D:Adm} we sketch some of the proofs:
\bp
Let $\xi\in \Adm (A)$ and $\eta\in  \Adm (B)$. Then there are rank one projections $P_j$ and $Q_j$ such that $A= \sum _j \xi_jP_j$ and $B= \sum _j\eta_jQ_j$. Then
\item [(vi)] 
$$A+B= \xi_1P_1+ \eta_1 Q_1+  \xi_2P_2+ \eta_2 Q_2+\cdots$$ Thus $\xi\oplus \eta\in \Adm(A+B)$. Considering $A=B$ we see that  $\Adm(A)\oplus \Adm(A)\subsetneq 2\Adm(A)$
\item [(viii)] $$A= tA+(1-t)A= t \xi_1P_1+ (1-t)\eta_1 Q_1+  t\xi_2P_2+ (1-t)\eta_2 Q_2+\cdots$$
 and thus $t\xi\oplus (1-t)\eta\in \Adm(A)$
\item [(ix)]
$$A= \eta_1P_1+ (\xi_1-\eta_1)P_1+  \eta_2P_2+ (\xi_2-\eta_2)P_2+\cdots+$$
and thus $\eta\oplus (\xi-\eta)\in \Adm(A).$
\ep 

Two positive operators $A,B\in  B(H)_+$ are called Murray-von Neumann equivalent if $B=VAV^*$ for some partial isometry $V$ with $V^*V\ge R_A$, (equivalently, if $B=XX^*$ and $A= X^*X$ for some $X\in B(H)$). Thus condition (ii) states that Murray-von Neumann equivalent operators have the same collection of admissible sequences. The converse holds if $A$ and $B$ are compact
\bR{R:MvN}
If $A,B\in  K(H)_+$, then $\Adm(A)=\Adm(B)$ if and only if $A$ and $B$ are Murray-von Neumann equivalent.
\eR
\bp
By the comments after Definition \ref {D:Adm} we have  $\{\xi\in \ell^\infty_+\mid \xi\preceq \lambda(A)\}=\Adm(A)=\Adm (B)=\{\xi\in \ell^\infty_+\mid \xi\prec \lambda(B)\} $. Thus $\lambda(A)\preceq \lambda (B)\preceq \lambda(A)$ and hence $\lambda(A)=\lambda(B)$. Thus $B=VAV^*$ for a partial isometry with $V^*V= R_A$ and hence $A$ and $B$ are Murray-von Neumann equivalent.
\ep
We can say more about admissible sequences of projections.

\bL{L:proj} Let $\xi\in \Adm (Q)$ for some projection $Q$ and let $P$ be a projection.
\item[(i)] $\xi\in \Adm (P)$ if and only if $\sum_j \xi_j= \tr(P).$  In particular, $\Adm(P)=\Adm(Q)$ if and only if $P$ and $Q$ are Murray-von Neumann equivalent.
\item[(ii)] $1-\xi\in \Adm (I-P)$ if and only if  $\sum_j(1-\xi_j) =\tr(I-P)$.  
\item[(iii)] $\diag\xi\in E(\mathscr U(P))$ if and only if $\sum_j \xi_j= \tr(P)$ and $\sum_j(1-\xi_j) =\tr(I-P)$.
\eL
\bp
Recall first that $\sum_{j=1}^\infty \xi_j= \tr(Q)$ by Lemma \ref {L: trivial properties} (vii) and that  $\diag \xi= E(R)$ for some projection $R\sim Q$. Then $\diag (1-\xi)= E(I-R)$ and in particular $1-\xi\in \Adm (I-R)$. The necessity part of (i) and (ii) follow from Lemma \ref {L: trivial properties} (vii), and the necessity of part (iii) is obvious.
\item[(i)] If $\sum_j \xi_j= \tr(P)$, then $\tr(P)= \tr(Q)$,   hence $P\sim Q$ and thus  $\Adm (P)=\Adm (Q)$ by Lemma \ref {L: trivial properties} (ii).
\item[(ii)] Since $1-\xi\in \Adm (I-R)$, by (i) it follows that $1-\xi\in \Adm (I-P)$.
\item [(iii)]  Since $\tr(P)=\sum_j \xi_j = \tr(R)$ and $\tr(I-P)= \sum_j(1-\xi_j) = \tr(I-R)$, it follows that $R$ is unitarily equivalent to $P$ and hence $\diag\xi\in E(\mathscr U(P))$.
\ep

It is worthwhile noticing that a sequence $\xi$ satisfies the Kadison condition if and only if $1-\xi$ satisfies the Kadison condition and that adding or deleting $0$ or $1$ entries from the sequence does not affect the Kadison condition. 
Notice also that if $\xi $ satisfies the Kadison condition, then either $\sum_j\xi_j=\infty$ or  $\sum_j\xi_j\in \mathbb N$. In either case there is always some projection $P$ such that $\tr(P)=\sum _j\xi_j$.

Before we present a proof of  Theorem \ref {T:main theorem}, we would like to explore a  consequence of it. Notice first that the Kadison condition is not necessary for a sequence to be admissible for  a sum of projections $A$ that is not a projection itself. Indeed  then $\|A\|>1$ and hence by Proposition \ref{P: bridge}, $\Adm(A)$ contains some sequences $\xi$ with values not bounded by 1.  Even sequences $\xi\in \Adm(A)$ that are bounded by 1 may fail to satisfy the Kadison condition as the following example illustrates. 
\begin{example} \label{E: not Kadison}
Let $A=2I+ \delta E_o$ with $E_o$ is a rank one projection and $0<\delta< 1$.  $A$ is a sum of projections because $\|A\|_{e}>1$ (see paragraph before Theorem \ref {T;sums of proj}). Since $2I$ is the sum of rank one projections,  $\xi:=\langle \delta, 1, 1, \cdots\rangle\in \Adm (A)$, but 
$ \sum\{\xi_i \mid \xi_i\le \frac{1}{2}\} = \delta$ and $ \sum\{(1-\xi_i) \mid \xi_i>   \frac{1}{2}\} = 0$
and thus 
$\xi$ does not satisfy the Kadison condition. 
\end{example}
The next lemma will help in finding a large class of admissible sequences of noncompact operators.
\bL{L:admiss by ess norm} Let $A\in B(H)_+\setminus K(H)$. Then every sequence $\xi$ such that  $0\le \xi_j \le \|A\|_{e}$ and $\xi_j < \|A\|$ for all $j$,  there are rank one projections $P_j$ such that $\sum_{j=1}^\infty \xi_jP_j\le A$. 
\eL

\bp
By scaling if necessary $A$, we can assume that $\|A\|_{e}=1$. If  $\chi_{[1,\infty)}(A)$ is infinite, then decompose it into an infinite sum of mutually orthogonal rank one projections $P_j$. Then 
$$A\ge A\chi_{[1,\infty)}(A)\ge \chi_{[1,\infty)}(A)=\sum_{j=1}^\infty P_j\ge \sum_{j=1}^\infty \xi_jP_j.$$
Thus assume that $\chi_{[1,\infty)}(A)$ is finite. Then $$\|A\chi_{(0,1)}(A)\|_e=\|A-A\chi_{[1,\infty)}(A) \|_e= \|A\|_e=1$$  and hence there is  a strictly increasing sequence $\alpha_n$ converging to $1$ and starting with $\alpha_1=0$, such that $\chi_{(\alpha_n,\alpha_{n+1})}(A)\ne 0$ for every $n$. Let $Q_n\le \chi_{(\alpha_n,\alpha_{n+1})}(A)$ be a rank one projection for every $n$. Then $$A= A\chi_{[1, \infty)}(A) + \sum_{n=1}^\infty A\chi_{(\alpha_n,\alpha_{n+1})}(A)\ge  A\chi_{[1, \infty)}(A)+\sum_{n=1}^\infty\alpha_nQ_n.$$
Next we consider two cases. If $\|A\|=1$ and hence $\xi_j<1$ for all $j$, then we can find a strictly increasing sequence of integers $n_j$ such that $\alpha_{n_j}\ge \xi_j$ for every $j$. Then $$A\ge \sum_{j=1}^\infty \xi_jQ_{n_j}.$$
If on the other hand $\|A\|>1$,  then the operator $k:= (A-I)\chi_{[1, \infty)}(A)$ is nonzero and has finite rank. Choose a  strictly increasing sequence of integers $n_j$ such that  $ \sum_{j=1}^\infty (1- \alpha_{n_j})\le \tr (k)$ and choose $0<t\le 1$ such that  $\tr(tk)=  \sum_{j=1}^\infty (1- \alpha_{n_j})$. Then
$$A\ge B:= \chi_{[1, \infty)}(A)+ tk +  \sum_{j=1}^\infty  \alpha_{n_j}Q_{n_j}.$$
Then $tk= (B-I)\chi_{[1, \infty)}(B)$, and $\sum_{j=1}^\infty (1- \alpha_{n_j})Q_{n_j}= (I-B)\chi_{[0, 1)}(B).$ Therefore by Theorem \ref {T;sums of proj}, $B=\sum_{j=1}^\infty P_j$ is a sum of infinitely many rank one projections $P_j$  and hence
$A\ge \sum_{j=1}^\infty \xi_jP_j$.
\ep
The following lemma shows that the  condition that $\xi_j\le\|A\|_e$ cannot be substantially weakened. 

\bL{L:ineq}
Let $0\le B\le A\in B(H)$. Then for every $\xi\in \Adm (B)$,  $\tr \big((A-I)_+\big) \ge \sum\{(\xi_j-1)\mid \xi_j>1\}$.
\eL
\bp
Let $\Lambda:=\{j\mid \xi_j>1\}$. Then 
$$A= (A-I)_+- (I-A)_+R_A + R_A\ge B= \sum_{j\in \Lambda}\xi_jP_j$$
for some rank-one projections $P_j$. 
Thus $C:=\sum_{j\in \Lambda}P_j$ converges strongly (if $\Lambda $ is infinite). To simplify notations , assume without loss of generality that $R_A=I.$ By Proposition \ref {P: bridge} there is a partial isometry $V$ with $V^*V= R_C$ and  such that $E(VC V^*)=I$. Then 
$$V(A-I)_+V^*- V(I-A)_+V^*+ VV^*\ge VCV^*+\sum_{j\in \Lambda}(\xi_j-1)VP_jV^*,$$
hence
$$
E(V(A-I)_+V^*)\ge E(VCV^*) -E(VV^*)+\sum_{j\in \Lambda}(\xi_j-1)E(VP_jV^*)\ge \sum_{j\in \Lambda}(\xi_j-1)E(VP_jV^*)
$$
whence the conclusion follows by computing the trace of both terms.
\ep
\bC{C: inside Adm}
Let $A\in B(H)_+\setminus K(H)$ and let $\xi$ be a sequence with $0\le \xi_j\le \|A\|_e$ and $\xi_j< \|A\|$ for all $j$ and such that $\sum \{\xi_j\mid \xi_j\le \alpha\}=\infty$ for some $0< \alpha< \|A\|_e$. Then $\xi\in \Adm(A)$.
\eC
\bp 
 Let $\Lambda:=\{j\mid \xi_j\le \alpha\}$. Decompose $A$ into a sum $A= A_1\oplus A_2$ with $\|A_1\|_{e}= \|A_2\|_{e}= \|A\|_{e}$ and $\|A_1\|= \|A\|$. Then by Lemma \ref {L:admiss by ess norm}, there are rank one projections $P_j$ such that $\sum_{j\not \in \Lambda}\xi_jP_j\le  A_1\le A.$  Let $B:=A-  \sum_{j\not \in \Lambda}\xi_jP_j$. Then $A_2\le B\le A$ and hence $\|B\|_e= \|A\|_e$. By Lemma \ref {L: trivial properties} (vi) it is enough to prove that $\langle\xi_j\rangle_{j\in \Lambda}\in\Adm (B)$.

 To simplify notations, assume that $\xi_j\le \alpha$ for every $j$, i.e., $\Lambda = \mathbb N$, (which corresponds to $B=A$). Let $\alpha< \beta < \|A\|_{e}$. Then $\|\frac{1}{\beta}A\|_{e}= \frac{\|A\|_{e}}{\beta}>1$ and hence $\frac{1}{\beta}A$ is a sum of projections by Theorem \ref {T;sums of proj}. Moreover,  $\tr \big(\frac{1}{\beta}A\big)=\sum_{j}\frac{\xi_j}{\beta}= \infty$. Then $\frac{\xi}{\beta}$ satisfies the conditions of Theorem \ref {T:main theorem} and hence $\frac{\xi}{\beta}\in \Adm \big(\frac{1}{\beta}A\big)$. But then $\xi\in \Adm(A)$ by Lemma \ref{L: trivial properties} (v).
 \ep

\section {The proof of the main theorem}

As we mentioned in the introduction, the proof of Theorem \ref {T:main theorem} is based for all but one case on the following result which by the transitivity of the $\prec$ relation is an obvious consequence of Theorem \ref {T:KkLd}.  

\bC{C:from KornLars} Let $\xi\in \mathbb R^n_+$ and $\eta\in \mathbb R^m_+$ with  $\xi\prec\eta$. Then for every collection of rank one projections $\{E_j\}_1^m$ there are rank one projections $\{P_j\}_1^n$ such that $\sum_{j=1}^n \xi_jP_j= \sum_{j=1}^m \eta_jE_j$.
\eC

To facilitate the application of this corollary, we present here two simple consequences of the definition of majorization.

\bL{L:elem} Let $0\le \xi_j\le 1$ for all $j$, and $\sum_{j=1}^\infty \xi= N+ r$ with $N\in \mathbb N\cup\{0\}$ and $0\le r< 1$. Then
\item [(i)] $\xi\prec\langle~ \overbrace{1,1, \cdots, 1}^N, r\rangle.$
\item [(ii)] If  $r=r_1+r_2$ with $0< r_2\le r_1$, then
$\xi\prec\langle~ \overbrace{1,1, \cdots, 1}^N, r_1,r_2\rangle$ if and only if 
$\sum_{j=1}^{N+1}\xi^*_j\le N+r_1$.
\eL
\bp
(i) Let $\eta= \langle \overbrace{1,1, \cdots, 1}^N, r, 0, \cdots \rangle$. Notice that $\eta=\eta^*$. For every integer $m\le N$,  $ \sum_{j=1}^m\xi^*_j\le m = \sum_{j=1}^m\eta_j$. For every $m>N$,   $ \sum_{j=1}^m\xi^*_j\le  \sum_{j=1}^\infty \xi^*_j
=N+r=  \sum_{j=1}^m\eta_j.$\\
(ii) $\eta= \langle \overbrace{1,1, \cdots, 1}^N, r_1,r_2, 0, \dots \rangle$. Here too $\eta=\eta^*$. As in (i), for every integer $m\le N$,  $ \sum_{j=1}^m\xi^*_j\le m = \sum_{j=1}^m\eta_j$ and for every  $m\ge N+2$, $$\sum_{j=1}^m\xi^*_j\le \sum_{j=1}^\infty\xi_j = N+r_1+r_2= \sum_{j=1}^{N+2}\eta_j= \sum_{j=1}^{m}\eta_j.$$
Thus for $\xi\prec \eta$ it is necessary and sufficient to have 
$\sum_{j=1}^{N+1}\xi^*_j\le\sum_{j=1}^{N+1}\eta_j= N+r_1.$
\ep
We first dispense of the finite rank case. Notice that if $0\le \xi_j\le  1$  for all $j$ and $\sum_{j=1}^\infty \xi_j\in \mathbb N$, then $\xi $ satisfies the Kadison condition.
\bL{L: fin rank}
Assume that $A= \sum_{j=1}^nE_j$ with $E_j$ rank one projections, $0\le \xi_j\le  1$ for all $j$, and $\sum_{j=1}^\infty \xi_j=\tr(A)$. Then $\xi\in \Adm (A)$
\eL
\bp
Let $m$ be an integer such that $\sum_{j=1}^m\xi_j= n-1+r$ with $0\le r< 1$. Then $\langle \xi_1, \cdots, \xi_m\rangle \prec \langle~ \overbrace{1,1, \cdots, 1}^{n-1}, r\rangle$ by Lemma \ref {L:elem}  (i) and hence by Corollary \ref {C:from KornLars}, there are rank one projections $P_j$ such that
$$\sum_{j=1}^{n-1}E_j+ rE_n= \sum_{j=1}^m\xi_jP_j$$ where we use the convention to drop a sum $\sum_{j=1}^{0}$. But then $1-r= \sum_{j=m+1}^\infty \xi_j$ and hence
$$\sum_{j=1}^m\xi_jP_j+ \sum_{j=m+1}^\infty \xi_jE_n= \sum_{j=1}^n E_j= A.$$
\ep 
From now on we assume that  $A= \sum_{j=1}^\infty E_j$ is an infinite sum of rank one projections, and that $\xi $ is a sequence satisfying the Kadison condition and such that $\sum_{j=1}^\infty \xi_j=\infty=\tr(A)$. For all but the last step in the proof we will further assume that assume that  $0< \xi_j <1$, and to  simplify notations, we decompose $\xi$ into two disjoint subsequences    $\xi= \mu\oplus (1-\lambda)$ where $\mu:=\langle\mu_j\rangle_1^M$ with $0\le \mu_i\le \frac{1}{2}$ and $\lambda:=\langle\lambda_j\rangle_1^N$,  with $0\le \lambda_i< \frac{1}{2}$, and $M, N\in \mathbb N \cup \{\infty\}$ ($0\le M, N\le \infty $ for short.) In terms of the sequence $\mu$ and $\lambda$, Kadison's  condition states that
either
$$ \sum_{j=1}^M \mu_j+ \sum_{j=1}^N \lambda_j = \infty\\
\quad \text{or}\quad 
 \sum_{j=1}^M \mu_j- \sum_{j=1}^N \lambda_j \in \mathbb Z.
$$ 
Notice that since $\xi= \mu\oplus (1-\lambda)$ satisfies the Kadison condition, then so does $1-\xi= \lambda\oplus (1-\mu)$, and that the roles of $M$ and $N$ are inverted.
Also,  
$$\sum_j\xi_j= \sum_{j=1}^M \mu_j+ \sum_{j=1}^N(1-\lambda_j )=\infty$$ hence if $\sum_{j=1}^M \mu_j< \infty$ then $N=\infty$.   

We will  prove the theorem by proving  the following complementary cases (for which we assume that $0< \xi_j< 1$ for all $j$):
\begin{itemize}
\item $\sum_{j=1}^M \mu_j = \infty$ (Lemma  \ref {L: infinite sum of mu}.)
\item $\sum_{j=1}^N \lambda_j = \infty$ (Lemma  \ref {L: infinite sum of lambda}.)\\
\noindent And the two possible cases when $\sum_{j=1}^M \mu_j < \infty$ and $\sum_{j=1}^N \lambda_j < \infty$:
\item $M=\infty$ and $N= \infty$ (Lemma \ref {L: finite with M=N infinite}.)
\item $M<\infty$ and $N= \infty$ (Lemma \ref  {L:M finite N infinite}.)
\end{itemize}
The proof of the theorem will then be completed by removing the condition that $0< \xi_j< 1$ for all $j$.

\bL{L: infinite sum of mu} If $\sum_{j=1}^\infty \mu_j=\infty$, then $\mu\oplus (1-\lambda) \in \Adm(A)$.
\eL
\bp
Assume that $N=\infty$ - the case when $N< \infty$ (including when $N=0$)  is similar and is left to the reader.

Since $\sum_{j=1}^\infty \mu_j=\infty$ there is some $n$, necessarily $n> 2$, such that $ \sum_{j=1}^{n} \mu_j\ge 1+\lambda_1$, and let  $n_1\in \mathbb N$ be the smallest such integer,  that is, $ \sum_{j=1}^{n_1} \mu_j +1-\lambda_1 = 2+ r_1\quad\text{for some }0\le r_1< \frac{1}{2}.$ 
By Lemma \ref {L:elem} (i) we see that
$\langle \mu_1, \cdots, \mu_{n_1}, 1-\lambda_1\rangle \prec \langle 1,1, r_1\rangle.$
By Corollary \ref {C:from KornLars}, there are rank one projections $Q_1$, and $P_j$ for $1\le j\le n_1$ such that
$A_1:= \sum_{j=1}^{n_1}\mu_j P_j + (1-\lambda_1)Q_1 = E_1+E_2 + r_1E_3.$

For the second step, let $n_2> n_1$, be such that
$ \sum_{j=n_1+1}^{n_2} \mu_j +1-\lambda_2= 2-r_1 +r_2$ for some $0\le r_2<\frac{1}{2}.$
We claim that
$\langle \mu_{n_1+1}, \cdots, \mu_{n_2}, 1-\lambda_2\rangle \prec \langle 1,1- r_1,r_2\rangle.$
Indeed, $n_2\ge n_1+2$ and if we let $\zeta$ be the monotone non increasing rearrangement of the finite sequence $\langle \mu_{n_1+1}, \cdots, \mu_{n_2}, 1-\lambda_2\rangle$, then $\zeta_1 = 1-\lambda_2< 1$ and $\zeta_1+\zeta_2= 1-\lambda_2+  \underset{n_1+1\le j\le n_2}{\max} \mu_j < \frac{3}{2}< 2-r_1.$
Thus the claim follows from Lemma \ref {L:elem} (ii).
Invoking again Corollary \ref {C:from KornLars}, we can find rank one projections $Q_2$, and $P_j$ for $n_1+1\le j\le n_2$ such that
$$A_2:=  \sum_{j=n_1+1}^{n_2}\mu_j P_j + (1-\lambda_2)Q_2 = (1-r_1)E_3+E_4+r_2E_5$$
and hence
$$A_1+A_2 =  \sum_{j=1}^{n_2}\mu_j P_j + \sum_{j=1}^2 (1-\lambda_j)Q_j= \sum_{j=1}^4E_j+ r_2E_5.$$
Iterating this construction we can find an increasing sequence of integers $n_k$, rank one projections $P_j$ and $Q_j$, and scalars $0\le r_k<\frac{1}{2}$ such that for every $k$ 
$$ \sum_{j=1}^{n_k}\mu_j P_j + \sum_{j=1}^k (1-\lambda_j)Q_j= \sum_{j=1}^{2k}E_j+ r_kE_{2k+1}.$$
By hypothesis, the series $\sum_{j=1}^n E_j\underset{s}{\to}A$, hence $E_j\underset{s}{\to} 0$ and thus  $$ A= \sum_{j=1}^{\infty}\mu_j P_j + \sum_{j=1}^\infty (1-\lambda_j)Q_j,$$
which concludes the proof.
\ep

\bL{L: infinite sum of lambda} If  $\sum_{j=1}^\infty \lambda_j=\infty$, then $\mu\oplus (1-\lambda) \in \Adm(A)$.
\eL
\bR{R: ovious1 if proj}
If $A=P$ is a projection, then Lemma \ref {L: infinite sum of lambda} is an immediate consequence of Lemma \ref {L: infinite sum of mu} and Lemma \ref {L:proj}.
\eR 
If $A$ is not a projection, then  we need first the following lemma.

\bL{L: infinite lambda}
If  $\sum_{j=1}^\infty \lambda_j=\infty$ and $B=\sum_{k=1}^Ks_kF_k$ where $K\in \mathbb N$, $0\le s_k< 1$, and $F_k$ are rank-one projections, then $1-\lambda\in\Adm(A+B)$.
\eL
\bp
Choose  $n_1$ such that $\sum_{j=1}^{n_1} \lambda_j \ge 2K+1$.  Since $\lambda_j<\frac{1}{2}$ for all $j$,
$$
\sum_{j=1}^{n_1}(1- \lambda_j)> \sum_{j=1}^{n_1} \lambda_j\ge K+1 +
 \sum_{k=1}^K s_k.$$ Then
\be{e:11} \sum_{j=1}^{n_1}(1- \lambda_j)=N_1+ \sum_{k=1}^K s_k+  r_1\ee
for an integer $N_1\ge K+1$ and $0\le r_1< 1$. Then
\be {e:n1}n_1=\sum_{j=1}^{n_1}(1- \lambda_j)+  \sum_{j=1}^{n_1} \lambda_j\ge  N_1 + \sum_{k=1}^K s_k+  r_1 + 2K+1.\ee
In particular, $n_1\ge N_1+K+1$. We claim that
$$\langle 1-\lambda_1, \cdots, 1-\lambda_{n_1}\rangle \,\prec \, \langle \overbrace{ 1, \dots,1 }^{N_1}, s_1, \cdots s_K, r_1\rangle.$$
Let $\eta: =  \langle \overbrace{ 1, \dots,1 }^{N_1}, s_1, \cdots s_K, r_1\rangle.$ For $1\le i\le N_1$, $\eta^*_i=\eta_i$ and hence
$\sum_{j=1}^i(1-\lambda_j)<  i= \sum_{j=1}^i\eta^*_j$.  For $i\ge N_1+K+1$ the majorization inequality follows from the equality in (\ref{e:11}). It remains to verify the majorization inequality for   ${N_1}< i\le N_1+K$. Then  
\begin{alignat*}{2}
\sum_{j=1}^i(1-\lambda_j)&= \sum_{j=1}^{n_1}(1-\lambda_j)- \sum_{j=i+1}^{n_1}(1-\lambda_j)\\
&< N_1 + \sum_{k=1}^Ks_k +r_1 - \frac{{n_1}-i}{2} &(\text{by (\ref {e:11}) and since $\lambda_j< \frac{1}{2}$)}\\
&< N_1 + \sum_{k=1}^Ks_k +r_1- \frac{1}{2}(N_1 + \sum_{k=1}^K s_k+  r_1 + 2K+1) + \frac{1}{2}( N_1+K)&(\text{by (\ref {e:n1})}\\
&=N_1 + \frac{1}{2}\big ( \sum_{k=1}^K s_k +r_1-K-1\big)\\
&< N_1&(\text{since  $s_k<1$ and $r_1< 1$})\\
&\le \sum_{j=1}^i \eta^*_j
\end{alignat*}
which completes the proof of the claim.
Then by Corollary \ref {C:from KornLars}, there are rank one projections $Q_j$ such that
$$\sum_{j=1}^{n_1}(1-\lambda_j)Q_j= \sum_{j=1}^{N_1}E_j+B+ r_1E_{{N_1}+1}.$$
Thus $A+B=\sum_{j=1}^{n_1}(1-\lambda_j)Q_j+ (1-r_1) E_{{N_1}+1}+ \sum_{j={N_1}+2}^\infty E_j.$

By applying the result just obtained to $(1-r_1) E_{{N_1}+1}+ \sum_{j={N_1}+2}^\infty E_j$,
we can find  $n_2>n_1$, $N_2>  N_1+1$,  rank one projections $Q_j$ and  and a remainder $0\le r_2<1$ such that 
$$\sum_{j=n_{1}+1}^{n_2}(1-\lambda_j)Q_j= (1-r_{1}) E_{{N_{1}}+1}+\sum_{j=N_{1}+2}^{N_2}E_j+ r_2E_{{N_2}+1}$$
and hence
$$
 \sum_{j=1}^{n_2}(1-\lambda_j)Q_j= B+ \sum_{j=1}^{N_2}E_j+ r_2E_{N_2+1}.
$$
Iterating, we find strictly increasing sequences of integers $n_k$ and $N_k$ and rank-one projections $Q_j$ such that
$$\sum_{j=1}^{n_k}(1-\lambda_j)Q_j=B+ \sum_{j=1}^{N_k}E_j+ r_kE_{{N_k}+1}.$$
Since $r_kE_{{N_k}+1}\underset{s}{\to} 0$ and $ \sum_{j=1}^{N_k}E_j\underset{s}{\to} A$ we obtain
$$\sum_{j=1}^{\infty}(1-\lambda_j)Q_j= A+B$$
which concludes the proof.

\ep
\noindent Now we can proceed with the proof of Lemma \ref {L: infinite sum of lambda}.
\bp
By Lemma \ref {L: infinite sum of mu}, we can assume without loss of generality that $\sum_{j=1}^\infty \mu_j< \infty$. 
Partition $\mathbb N$ in a disjoint union of $K$ sets $\Lambda_k$ such that $\sum_{j\in \Lambda_k}\mu_j\le 1$. Set $s_k:=1- \sum_{j\in \Lambda_k}\mu_j$ and $B:= \sum_{k=1}^Ks_kE_k$. Then by Lemma \ref {L: infinite lambda} there are rank one projections $Q_j$ such that
$$\sum_{j=1}^\infty (1-\lambda_j)Q_j =  \sum_{k=1}^Ks_kE_k+  \sum_{j=K+1}^\infty E_k.$$
But then
$$A=   \sum_{k=1}^\infty E_k= \sum_{k=1}^K \sum_{j\in \Lambda_k}\mu_jE_k+ \sum_{j=1}^\infty (1-\lambda_j)Q_j,$$
which completes the proof.
\ep
Thus we are reduced to consider the case where $\sum_{j=1}^N \lambda_j+\sum_{j=1}^M \mu_j<\infty$ but $\sum_{j=1}^N (1-\lambda_j)=\infty$ and hence $N=\infty$.

\bL{L: finite with M=N infinite} Assume that $M=N=\infty$ and $\sum_{j=1}^\infty \lambda_j<\infty$, $\sum_{j=1}^\infty \mu_j<\infty$.  Then  $\mu\oplus (1-\lambda) \in \Adm(A)$.
\eL
\bp
 Let $k:= \sum_{j=1}^\infty \lambda_j-\sum_{j=1}^\infty \mu_j$. Then $k\in \mathbb Z$ by hypothesis. Choose a positive integer $n_1\ge  k+1$ such that $ \sum_{j=n_1+1}^\infty \lambda_j<\frac{1}{2}$.  Since $\lambda_j\ne 0$ for all $j$, choose $m_1$ such that $\sum_{j=m_1+1}^\infty \mu_j\le \sum_{j=n_1+1}^\infty \lambda_j$. Set $r_1:= \sum_{j=n_1+1}^\infty \lambda_j-\sum_{j=m_1+1}^\infty \mu_j$. Then $0\le r_1< \frac{1}{2}$ and 
$$\sum_{j=1}^{m_1}\mu_j+ \sum_{j=1}^{n_1}(1-\lambda_j)= n_1-k+r_1.$$
We have by Lemma \ref {L:elem} (i) that 
$$\langle\mu_1, \cdots, \mu_{m_1}, 1-\lambda_1, \cdots, 1-\lambda_{n_1}\rangle\,\prec \, \langle\overbrace{ 1, \dots,1 }^{n_1-k}, r_1\rangle.$$
Hence there are rank-one projections $P_j$ and $Q_j$ such that
$$
\sum_{j=1}^{m_1}\mu_jP_j+ \sum_{j=1}^{n_1}(1-\lambda_j)Q_j= \sum_{j=1}^{n_1-k}E_j+ r_1E_{n_1-k+1}.$$
For the next step,  choose $n_2>n_1+1$ and $m_2>m_1$ such that 
$$\sum_{j=m_2+1}^\infty \mu_j\le\sum_{j=n_2+1}^\infty \lambda_j\le \sum_{j=m_1+1}^\infty \mu_j.$$ Set $r_2= \sum_{j=n_2+1}^\infty \lambda_j- \sum_{j=m_2+1}^\infty \mu_j$.
Then $0\le r_2< \frac{1}{2}$ and
$$\sum_{j=m_1+1}^{m_2}\mu_j+ \sum_{j=n_1+1}^{n_2}(1-\lambda_j) = n_2-n_1 +r_2-r_1.$$
We claim that
$$\langle \mu_{n_1+1}, \cdots , \mu_{n_2}, 1- \lambda_{n_1+1}, \cdots, 1- \lambda_{n_2}\rangle\,\prec\, \langle\overbrace{ 1, \dots,1 }^{n_2-n_1-1}, 1-r_1, r_2\rangle.$$
Indeed, this follows from Lemma \ref {L:elem} (ii)  since $1-\lambda_j>\mu_i$ for every $i$ and $j$ and
\ba \sum_{j=n_1+1}^{n_2} 1- \lambda_j& = n_2-n_1- \sum_{j=n_1+1}^\infty \lambda_j + \sum_{j=n_2+1}^\infty \lambda_j \\
&\le n_2-n_1- \sum_{j=n_1+1}^\infty \lambda_j+\sum_{j=m_1+1}^\infty \mu_j\\&= n_2-n_1-r_1.
\end{align*}

Then by Corollary \ref {C:from KornLars}  there are rank one projections $P_j$ and $Q_j$ such that
$$
\sum_{j=m_1+1}^{m_2}\mu_jP_j+\sum_{j=n_1+1}^{n_2} (1- \lambda_j)Q_j= (1-r_1)E_{n_1-k+1}+ \sum_{j= n_1-k+2}^{n_2-k}E_j+ r_2E_ {n_2-k+1}$$
and hence
$$
\sum_{j=1}^{m_2}\mu_jP_j+\sum_{j=1}^{n_2} (1- \lambda_j)Q_j=  \sum_{j= 1}^{n_2-k}E_j+ r_2E_ {n_2-k+1}.
$$
Iterating this construction we get 
$$
\sum_{j=1}^{\infty}\mu_jP_j+\sum_{j=1}^{\infty} (1- \lambda_j)Q_j=  \sum_{j= 1}^{\infty}E_j=A.$$
\ep

The last case is 
\bL{L:M finite N infinite} Assume that $M<\infty$, $N= \infty$, and $\sum_{j=1}^\infty \lambda_j<\infty$.  Then  $\mu\oplus (1-\lambda) \in \Adm(A)$.
\eL
\bR{R: ovious2 if proj}
If $A=P$ is a projection, then Lemma \ref {L:M finite N infinite} is an immediate consequence of  
Lemma \ref {L:proj}. Indeed $1-\xi= \lambda\oplus (1-\mu)$ satisfies the Kadison condition and $\sum_{j=1}^\infty(1-\xi_j)< \infty$,  thus $\sum_{j=1}^\infty(1-\xi_j)=n$ for some $n\in \mathbb N$. By Lemma \ref {L: fin rank},  $1-\xi\in \Adm(Q)$ for any projection $Q$ with rank $n$, and hence by Lemma \ref {L:proj}, $\xi\in \Adm(P)$.
 \eR
If $A$ is not a projection, then Lemma \ref {L:proj} cannot be invoked and thus  a proof based on Corollary \ref {C:from KornLars} does not seem to be available. We will instead need  the estimate provided in the following  $2 \times 2$ matrix decomposition.

\bL {L: 2 x 2} Let $\eta, \xi\in \mathbb R_+^2$ with $\eta_1\ne \eta_2$ and $ \xi\prec \eta$, and let $u, u'\in \mathbb C^2$  be unit vectors. Then there are unit vectors $w, w'\in \mathbb C^2$ with   $w= \sigma u+ \tau u'$ for some scalars $\sigma, \tau \in \mathbb C$ such that
\item [(i)] $\xi_1w\otimes w+ \xi_2w'\otimes w'=\eta_1u\otimes u+ \eta_2u'\otimes u'$;
\item [(ii)] $|\sigma|^2+|\tau|^2\le 1$
\item [(iii)] $|\sigma|^2\le   \frac{\eta_1(\eta_1-\xi_2)}{\xi_1(\eta_1-\eta_2)}.$
\eL
\bp
If $\xi_1=\eta_1$ then $\xi_2=\eta_2$ and hence we can choose $w=u$ (and hence $\sigma=1$, $\tau=0$) and $w'=u'$. Then (i) holds trivially and we have equality in (ii)  and (iii).\\ If $\xi_2=\eta_1$ then $\xi_1=\eta_2$ and hence we can choose $w=u'$ (and hence $\sigma=0$, $\tau=1$) and $w'=u$.  Again, (i) holds trivially and we have equality in (ii)  and (iii).
Thus assume henceforth that $\xi_1\ne \eta_1$ and  $\xi_2\ne \eta_1$. Equally trivial is the case when $u\otimes u=u'\otimes u'$, thus let $\gamma: =|(u,u')|$ and assume that
$0\le \gamma <1$. We set $w'= \sigma' u+ \tau' u'$  and we claim that we can satisfy conditions (i)-(iii) while further assuming that   $\sigma>0$, $\sigma'>0$, $\tau>0$, and $\tau'<0$.

For the vector $w$ and $w'$ to have unit norm, it is sufficient  (and necessary)  that 
\begin{align}\sigma^2\tau^2 + 2 \gamma \sigma \tau = 1\label{e: 1st ellipse in x y}\\
\sigma'^2+\tau '^2 + 2 \gamma \sigma'\tau ' = 1\label{e: ellipse in x'y'}.
\end{align}
For (i) to hold, it  is sufficient (and necessary) 
to have  both
\ba
 \eta_1 u&+ \gamma \eta_2 u'\\
 &=\big(\eta_1 u\otimes u+ \eta_2 u'\otimes u' \big)u\\
 &= \big(\xi_1 w\otimes w+ \xi_2 w'\otimes w' \big)u\\
 &=\xi_1(\sigma + \gamma \tau)w+ \xi_2(\sigma' + \gamma \tau')w'\\
&=\Big( \xi_1 (\sigma^2+\gamma \sigma\tau )+  \xi_2 (\sigma'^2+\gamma \sigma '\tau ')\Big)u+ \Big( \xi_1 (\sigma \tau +\gamma \tau ^2)+  \xi_2 (\sigma '\tau '+\gamma \tau '^2)\Big)u'\\
\intertext{and}
 \gamma\eta_1 &u+ \eta_2 v\\
&=\big(\eta_1 u\otimes u+ \eta_2 u'\otimes u' \big)u'\\
 &= \big(\xi_1 w\otimes w+ \xi_2 w'\otimes w' \big)u'\\
&= \xi_1(\gamma \sigma + \tau)w+ \xi_2(\gamma \sigma' + \tau')w'\\
&=\Big( \xi_1 (\gamma \sigma ^2+\sigma \tau )+  \xi_2 (\gamma \sigma '^2+ \sigma '\tau ')\Big)u+ \Big( \xi_1 (\gamma \sigma \tau + \tau ^2)+  \xi_2 (\gamma \sigma '\tau '+\tau '^2)\Big)u'.
\end{align*}
Since $u$ and $v$ are linearly independent, (ii) holds if and only if the following system of four  equations is satisfied:
\begin{align*}
 \xi_1 (\sigma ^2+\gamma \sigma \tau )+  \xi_2 (\sigma '^2+\gamma \sigma '\tau ')&=\eta_1  \\
  \xi_1 (\sigma \tau +\gamma \tau ^2)+  \xi_2 (\sigma '\tau '+\gamma \tau '^2)&=\gamma \eta_2\\
   \xi_1 (\gamma \sigma ^2+\sigma \tau )+  \xi_2 (\gamma \sigma '^2+ \sigma '\tau ')&= \gamma \eta_1\\
   \xi_1 (\gamma \sigma \tau + \tau ^2)+  \xi_2 (\gamma \sigma '\tau '+\tau '^2)&= \eta_2.
\end{align*}
Since $\gamma\ne 1$, easy algebraic manipulations show that the above system  and hence (i) are equivalent to the system 
\be{e:second system} \begin{cases}
\xi_1\sigma \tau +\xi_2 \sigma '\tau '&=0\\
 \xi_1 \sigma ^2+  \xi_2 \sigma '^2&=\eta_1\\
   \xi_1 \tau ^2+  \xi_2 \tau '^2&= \eta_2.
\end{cases}
\ee
If we choose $\sigma '^2$ and $\tau '^2$ so to satisfy the second and third equation, then $$\xi_2^2\sigma '^2\tau '^2= (\eta_1- \xi_1\sigma ^2)(\eta_2- \xi_1\tau ^2)=\xi_1^2\sigma ^2\tau ^2 + \eta_1\eta_2- \xi_1( \eta_2\sigma ^2 + \eta_1\tau ^2).$$  Given that $\sigma, \tau, \sigma'>0$ and $\tau'< 0$, the system (\ref {e:second system}) is thus equivalent to
\be{e:third system} \begin{cases}
\eta_2\sigma ^2 + \eta_1\tau ^2&=  \frac{\eta_1\eta_2}{\xi_1}\\
 \xi_1 \sigma ^2+  \xi_2 \sigma '^2&=\eta_1\\
   \xi_1 \tau ^2+  \xi_2 \tau '^2&= \eta_2.
\end{cases}
\ee
Furthermore, if (\ref {e: 1st ellipse in x y}) and (\ref {e:third system}) hold   then also (\ref {e: ellipse in x'y'}) holds. indeed, 
 \ba \sigma '^2+\tau '^2 + 2 \gamma \sigma '\tau '&= \frac{\eta_1}{\xi_2}-  \frac{\xi_1}{\xi_2} \sigma ^2+ \frac{\eta_2}{\xi_2}-  \frac{\xi_1}{\xi_2} \tau ^2 - 2\gamma\frac{\xi_1}{\xi_2}\sigma \tau \\
&=\frac{\eta_1+\eta_2}{\xi_2}-\frac{\xi_1}{\xi_2}\big( \sigma ^2+\tau ^2 + 2 \gamma \sigma \tau \big)\\
&= \frac{\xi_1+\xi_2}{\xi_2}-\frac{\xi_1}{\xi_2} =1
\end{align*}

Thus if $\sigma$ and $\tau$ satisfy the system
\be{e:intersection} \begin{cases} 
\sigma ^2+\tau ^2 + 2 \gamma \sigma \tau  &= 1\\
 \eta_2\sigma ^2 + \eta_1\tau ^2&=  \frac{\eta_1\eta_2}{\xi_1}.
\end{cases}
\ee
and $\sigma '$ and $\tau '$ are chosen as above, then $w$ and $w'$ are unit vectors satisfying (i). We are going to prove that the (unique) solution of the system (\ref{e:intersection}) with $\sigma>0$ and $\tau>0$ satisfies also (ii) and (iii). Squaring within the first equation and  eliminating $\tau ^2$ we obtain
$$
4\gamma^2\sigma ^2( \frac{\eta_2}{\xi_1}- \frac{\eta_2}{\eta_1}\sigma ^2)= \big ( 1-\sigma ^2-\frac{\eta_2}{\xi_1}+ \frac{\eta_2}{\eta_1}\sigma ^2 \big)^2. 
$$
and hence 
\be{e:quadr raw}
4\gamma^2 \frac{\eta_2}{\eta_1}\sigma ^2\big( \frac{\eta_1}{\xi_1}-\sigma ^2\big) = \frac{(\eta_1-\eta_2)^2}{\eta_1^2}\Big( \frac{\eta_1(\xi_1-\eta_2)}{\xi_1(\eta_1-\eta_2)}-\sigma ^2\Big)^2.
\ee
Set \begin{align*} z&:=\sigma ^2\\
z_o&:= \frac{\eta_1(\eta_1-\xi_2)}{\xi_1(\eta_1-\eta_2)}\\
h&:= \frac{4\eta_1\eta_2\gamma^2}{(\eta_1-\eta_2)^2}\\
\alpha&:=  \frac{\eta_1-\eta_2}{\eta_1-\xi_2},\end{align*} and using the fact that $\eta_1+\eta_2= \xi_1+\xi_2$, we can rewrite (\ref {e:quadr raw}) as 
\be{e:quadratic} (1+h)z^2- (2+\alpha h)z_oz+z_o^2=0.\ee
Since $h\ge 0$, $\alpha >1$, and $z_o>0$, this equation has two real (and positive roots). Let $z_-$ be the smallest root, i.e., $$z_-= \frac{2z_o}{2+\alpha h+ \sqrt{ 4(\alpha-1)h+ \alpha^2h^2}}$$
Then $z\le z_o$.
It is now easy to verify  that $\sigma = \sqrt {z_-}$ and $\tau = \sqrt{\frac {\eta_2 }{\xi_1}- \frac {\eta_2 }{\eta_1}z_-}$ are indeed solutions of the system (\ref {e:intersection}) and that they satisfy (ii) and (iii). 

\ep

\bL{L:key case}
Assume that $\sum_{j=1}^\infty\lambda_j< 1$ and let $$B:= (1- \sum_{j=1}^\infty\lambda_j) u_1\otimes u_1 + \sum_{j=2}^\infty u_j\otimes u_j$$  for some unit vectors $u_j$ for which the series $\sum_{j=2}^\infty u_j\otimes u_j$ converges in the SOT. Then $1-\lambda\in \Adm(B)$.
\eL
\bp
Apply Lemma \ref  {L: 2 x 2} to  $B_2:= (1- \sum_{j=1}^\infty\lambda_j) u_1\otimes u_1 +   u_2\otimes u_2$  where $u=u_1$, $\eta_1=1- \sum_{j=1}^\infty\lambda_j$, $u'=u_2$, $\eta_2=1$,  $\xi_1= 1- \sum_{j=2}^\infty\lambda_j$, and $\xi_2= 1-\lambda_1$. Then there are unit vectors $w_2$(=$w$) and $v_1$(=$w'$) such that
$$B_2=  (1-\lambda_1)v_1\otimes v_1+ (1- \sum_{j=2}^\infty\lambda_j) w_2\otimes w_2 $$ where $$w_2= \sigma_2 u_1+ \tau_2 u_2,$$  
$|\sigma_2|^2+ |\tau_2|^2\le 1$, and $$|\sigma_2|^2\le \frac{\big(1- \sum_{j=1}^\infty\lambda_j\big) \big(\sum_{j=2}^\infty\lambda_j\big)}{\big(\sum_{j=1}^\infty\lambda_j\big)\big(1- \sum_{j=2}^\infty\lambda_j\big)}.$$
Next, apply Lemma  \ref  {L: 2 x 2} to  $B_3:= (1- \sum_{j=2}^\infty\lambda_j) w_2\otimes w_2+u_3\otimes u_3$. Then again there are unit vectors $w_3$ and $v_2$ such that
$$B_3=  (1-\lambda_2)v_2\otimes v_2+ (1- \sum_{j=3}^\infty\lambda_j) w_3\otimes w_3 $$ where $$w_3= \sigma_3 w_2+ \tau_3 u_3= \sigma_3\sigma_2u_1+ \sigma_3\tau_2u_2+ \tau_3u_3,$$ 
$|\sigma_3|^2+ |\tau_3|^2\le 1$, and $$|\sigma_3|^2\le \frac{\big(1- \sum_{j=2}^\infty\lambda_j\big) \sum_{j=3}^\infty\lambda_j}{\sum_{j=2}^\infty\lambda_j\big(1- \sum_{j=3}^\infty\lambda_j\big)}.$$
Therefore
$$(1- \sum_{j=1}^\infty\lambda_j) u_1\otimes u_1 +   u_2\otimes u_2+u_3\otimes u_3= (1-\lambda_1)v_1\otimes v_1+(1-\lambda_2)v_2\otimes v_2 +(1- \sum_{j=3}^\infty\lambda_j) w_3\otimes w_3.$$
Iterating the construction, we find a sequence of unit vectors $v_j$ and $w_j$ such that for every $n$
$$(1- \sum_{j=1}^\infty\lambda_j) u_1\otimes u_1+ \sum_{j=2}^{n}u_j\otimes u_j=  \sum_{j=1}^{n-1}(1-\lambda_j)v_j\otimes v_j+ \big(1- \sum_{j=n}^\infty \lambda_j \big)w_n\otimes w_n,$$
where \be{e: recurr} w_n= \sigma_n w_{n-1}+ \tau_n u_n,\ee $|\sigma_n|^2+ |\tau_n|^2\le 1$, and $$|\sigma_n|^2\le \frac{\big(1- \sum_{j=n-1}^\infty\lambda_j\big) \sum_{j=n}^\infty\lambda_j}{\sum_{j=n-1}^\infty\lambda_j\big(1- \sum_{j=n}^\infty\lambda_j\big)}.$$
Since by hypothesis the series  $\sum_{j=2}^{\infty}u_j\otimes u_j$ converges in the strong operator topology, to show that the series  $\sum_{j=1}^{\infty}(1-\lambda_j)v_j\otimes v_j$ converges to $B$, it is enough to show that $\big(1- \sum_{j=n}^\infty \lambda_j \big)w_n\otimes w_n\to 0$. As $1- \sum_{j=n}^\infty \lambda_j\to 1$, we need to prove that $w_n\underset{w}{\to} 0 $.
Solving the recurrence (\ref {e: recurr}) yields
$$w_n= \prod_{i=2}^n\sigma_i u_1+ \tau_2\prod_{i=3}^n\sigma_i u_2+\cdots+  \tau_k\prod_{i=k+1}^n\sigma_i u_k+\cdots+\tau_{n-1}\sigma_nu_{n-1}+ \tau_nu_n.$$
Let $\{e_n\}$ be an orthonormal basis of the Hilbert space and let 
$$x_n:=  \prod_{i=2}^n\sigma_i e_1+ \tau_2\prod_{i=3}^n\sigma_i e_2+\cdots+  \tau_k\prod_{i=k+1}^n\sigma_i e_k+\cdots+\tau_{n-1}\sigma_ne_{n-1}+ \tau_ne_n.$$
The sequence $\{x_n\}$ is bounded, indeed \ba\|x_n\|^2&=\prod_{j=2}^n|\sigma_j |^2+ |\tau_2|^2\prod_{j=3}^n|\sigma_j|^2+\cdots+  |\tau_k|^2\prod_{j=k+1}^n|\sigma_j |^2+\cdots+|\tau_{n-1}|^2|\sigma_n|^2+ |\tau_n|^2\\
&\le \prod_{j=2}^n|\sigma_j |^2+ (1-|\sigma_2|^2)\prod_{j=3}^n|\sigma_j|^2+\cdots+  (1-|\sigma_k|^2)\prod_{j=k+1}^n|\sigma_j |^2+\cdots+1- |\sigma_n|^2\\
&=  \prod_{j=2}^n|\sigma_j |^2+  \prod_{j=3}^n|\sigma_j |^2-  \prod_{j=2}^n|\sigma_j |^2 +\cdots+  \prod_{j=k+1}^n|\sigma_j |^2- \prod_{j=k}^n|\sigma_j |^2+\cdots+ |\sigma_n|^2-  \prod_{j=n-1}^n|\sigma_j |^2+1- |\sigma_n|^2\\
&=1.
\end{align*}

Moreover, for every $k$, 
\ba
|(x_n, e_k)|^2&= |\tau_k|^2\prod_{i=k+1}^n|\sigma_i|^2
\le \prod_{j=k+1}^n\frac{\big(1- \sum_{j=i-1}^\infty\lambda_j\big)\big( \sum_{j=i}^\infty\lambda_j\big)}{\big(\sum_{j=i-1}^\infty\lambda_j\big)\big(1- \sum_{j=i}^\infty\lambda_j\big)}\\
&= \frac{\big(1- \sum_{j=k}^\infty\lambda_j\big) \big(\sum_{j=n}^\infty\lambda_j\big)}{\big(\sum_{j=k}^\infty\lambda_j\big)\big(1- \sum_{j=n}^\infty\lambda_j\big)} \to 0
\end{align*}
Therefore, $x_n\underset{w}{\to} 0 $. Since $w_n= T^*x_n$ where $T^*:= \sum_{j=1}^\infty u_j\otimes e_j$ is bounded (in fact, $T^*$ is the synthesis operator of the Bessel sequence $\{u_j\}$,) it follows that $w_n\underset{w}{\to} 0 $, which concludes the proof.
\ep
We can now provide the proof of Lemma \ref {L:M finite N infinite} 
\bp

Let $k= \sum_{j=1}^\infty \lambda_j- \sum _{j=1}^M \mu_j$ and choose $n> k+1$ such that $r:=\sum_{j=n+1}^\infty\lambda_j< 1$.
Then
$$\sum_{j=1}^M\mu_j+\sum_{j=1}^n(1-\lambda_j)= n-k-r$$
and hence 
$$\langle \mu_{1}, \cdots , \mu_{m}, 1- \lambda_{1}, \cdots, 1- \lambda_{n}\rangle\,\prec\, \langle\overbrace{ 1, \dots,1 }^{n-k-1}, 1-r\rangle.$$
Then there are rank one projections $P_j$ and $Q_j$
$$
 \sum_{j=1}^M\mu_jP_j+ \sum_{j=1}^{n} (1- \lambda_j)Q_j= \sum_{j=1}^{n-k}E_j- rE_{n-k}.$$
By Lemma \ref {L:key case}, there are also rank-one projections $Q_j$ such that
$$\sum_{j=n+1}^\infty (1-\lambda_j)Q_j=rE_{n-k}+ \sum_{j=n-k}^\infty E_j$$
Combining these two decompositions we obtain
$$ \sum_{j=1}^M\mu_jP_j+ \sum_{j=1}^{\infty} (1- \lambda_j)Q_j=A$$
which completes the proof.
\ep

\bp [End of the proof of Theorem \ref {T:main theorem}]

Assume there are $0\le m\le \infty$ (resp., $0\le n\le \infty$)  indices for which $\xi_j=0$ (resp., $\xi_j=1$). Then $\xi= \tilde{\xi}\oplus 0_m\oplus 1_n$ where $\tilde{\xi}$ is the sequence obtained by dropping from $\xi$ all the entries that are $0$ or $1$,  $0_m:= \{0\}_1^m$, and  $1_n:=\{0\}_1^n$, and we adopt the convention of dropping a direct summand if the sequence is empty (has no entries).
We leave to the reader the trivial case where  $\tilde{\xi}$ or $1_n$ are empty. If $\xi$ satisfies the Kadison condition, then so does $\tilde{\xi}$. Let $t= \sum\{\xi_j\mid \xi_j< 1\}$, i.e., $t= \sum \tilde{\xi}_j$,  then $\sum \xi_j= t+n$ and $t\in \mathbb N\cup \{\infty\}$. For every sum of  projections $A$ such that $\tr(A)= t+n $, i.e., a sum of $t+n\le \infty$ rank-one projection, we can decompose $ A= A_t+A_n$ into the sum  $A_t$  of $t$ rank-one projections and $A_n$ the sum of $n$ rank-one projections. By the previous result, $\tilde{\xi}\in \Adm(A_t)$ and trivially, $1_n\in \Adm(A_n)$. Thus by Lemma \ref {L: trivial properties}, $\xi \in \Adm (A).$ 

\ep


\begin{thebibliography}{99}
\bibitem{AMRS}
Antezana,  J.,  Massey, P.,  Ruiz, M., and Stojanoff, D., \textit{The Schur--Horn Theorem for operators and frames with prescribed norms and frame operator,} Illinois J of Math.  to appear.

\bibitem{Am15} Argerami, M., {\it Majorisation and the Carpenter Theorem}, Int. Eq. Oper Theory, ({\bf 82}), 1, (2015), 33-49.

\bibitem{AW07} Arveson, W.,  {\it Diagonals of normal operators with spectrum}, PNAS {\bf 104.4} (2007), 1152-1158.
 
\bibitem{AK02}
Arveson, W.  and Kadison, R.V.,  \textit{Diagonals of self-adjoint operators.}  Operator theory, operator algebras, and applications, 
Contemp. Math. \textbf{414}  Amer. Math. Soc., Providence, RI (2006), 247--263.


\bibitem{BJdiag}Bownik, M.  and Jasper, J., \textit{ Diagonals of self-adjoint operators with finite spectrum}, Bull. Polish Acad. Sci. \textbf{63}, (2015), 249-260.

\bibitem{BJSchHor}Bownik, M.  and Jasper, J., \textit{ The Schur--Horn Theorem for operators with finite spectrum}, Trans. AMS, \textbf{367}, (2015), 5099-5140.

\bibitem{BJCarpenter}Bownik, M.  and Jasper, J., \textit{Constructive proof of Carpenter's theorem}, Canadian Math. Bull. ({\bf 57}), 3, 463-476.


\bibitem{CFKLT06} Casazza, P., Fickus M., Kovacevic, J., Leon, M., and Tremain, J., {\it A physical interpretation of tight frames,} Harmonic analysis and applications, Appl. Numer. Harmon. Anal., Birkhauser Boston, Boston, MA (2006), 51-76

\bibitem{CL02}
Casazza, P. and Leon, M., \textit{Frames with a given frame operator} (2002) Preprint.

\bibitem{CL10} Casazza, P. and Leon, M., \textit{Existence and construction of finite frames with a given frame operator} Int. J. Pure Appl Math, {\bf 63}, 2 (2010), 149-157.

\bibitem{ChoiWu2014} Choi, M.D.  and  Wu, P.Y., \textit{Sums of orthogonal projections}, J. Func. Anal., 267 (2014), 384-404

\bibitem{DraganKaftal2016} Dragan, C. and Kaftal, V.,  \textit{Sums of equivalent sequences of positive operators in von Neumann factors}, Houston Journal of Mathematics, {\bf42}, (2016), No. 3, 991-1017.

 \bibitem {DFKLOW}
 Dykema, K., Freeman, D., Kornelson, K., Larson, D., Ordower, M., and Weber, E.,
 \textit{Ellipsoidal tight frames and projection decompositions of operators} Illinois J of Math. \textbf{48}, (2004), 477-489.


\bibitem{Fp69}
Fillmore, P.\textit{ On sums of projections.} J. Funct.
Anal.\textbf{4} (1969), 146-152.

\bibitem{GiMa64}
Gohberg, I. C. and Markus, A. S., \textit{Some relations between eigenvalues and matrix elements of linear operators} Mat. Sb. \textbf{64} (106) (1964), 481-496 (Russian); Amer. Math. Soc. Transl. (2) \textbf{52} (1966) 201-216  (English)

\bibitem{Ha54}
Horn, R. A. \textit{Doubly stochastic matrices and the diagonal of a rotation matrix.} Amer. J. Math.  \textbf{76} (1954), 620-630.

\bibitem{Kr02a}
Kadison, R.,  \textit{The Pythagorean Theorem I: the finite case,}  Proc. Natl. Acad. Sci. USA \textbf{99} (7) (2002), 4178-4184.

\bibitem{Kr02b}
Kadison, R.,  \textit{The Pythagorean Theorem II: the infinite discrete case,}  Proc. Natl. Acad. Sci. USA \textbf{99}  8  (2002), 5217-5222.

\bibitem{KvLj08} Kaftal, V. and Loreau, J., {\it Kadison's Pythagorean Theorem and essential codimension}, Int Eq. Oper. Th (to appear)


\bibitem{KNZ strong sums}
Kaftal, V., Ng, P.W., and  Zhang, S. \textit{Strong sums of projections in von Neumann factors,} J. 
Funct.  Anal. {\bf 257} (2009) 8, 2497-2529


\bibitem{KNZFiniteSumsVNA} Kaftal, V., Ng, P.W., and  Zhang, S. \textit{Finite sums of projections in von Neumann algebras,}
Trans. Amer. Math. Soc., \textbf{365} (2013) 5, 2409-2445.


\bibitem{KaftalWeiss2010}  Kaftal, V.  and Weiss, G.,
{\it An infinite dimensional Schur--Horn Theorem and majorization theory.}
J. Funct. Anal. \textbf{259} (2010), 12, 3115-3162.


\bibitem{KkLd04}
Kornelson, K. and Larson, D., \textit{Rank-one decompositions of operators and construction of frames,}
Wavelets, frames and operator theory, Contemp. Math \textbf{345}, Amer. Math. Soc. Providence, Ri (2004),  203-214.  

\bibitem{LjWg15}Loreaux, J. and  Weiss, G., {\it Majorization and a Schur--Horn Theorem for positive compact operators, the nonzero kernel case}, J. Func Anal, {\bf 268}, (2015), 3, 703-731. 


\bibitem{aM64}
Markus, A.S., \textit{The eigen-and singular values of the sum and product of linear operators,} Uspekhi Mat. Nauk 4  \textbf{118} (1964),  93-123.





\bibitem{Si23}
Schur, I., \textit{\"{U}ber eine Klasse von Mittlebildungen mit Anwendungen auf der Determinantentheorie,} Sitzungsber. Berliner Mat. Ges., \textbf{22}  (1923),  9-29.

\end{thebibliography}
 \end{document}